\documentclass{amsart}

\usepackage[latin1]{inputenc}
\usepackage[T1]{fontenc}
\usepackage{enumerate}
\usepackage{diagbox}  

\usepackage[english]{babel}

\usepackage{amsmath,amsfonts,amssymb,amsthm,amscd}

\newtheorem{theorem}{Theorem}[section]
\newtheorem{lemma}[theorem]{Lemma}
\newtheorem{corollary}[theorem]{Corollary}
\newtheorem{proposition}[theorem]{Proposition}
\newtheorem{claim}{Claim}

\theoremstyle{definition}
\newtheorem{definition}[theorem]{Definition}

\newtheorem{problem}{Problem}

\theoremstyle{remark}
\newtheorem{remark}[theorem]{Remark}

\numberwithin{equation}{section}

\newcommand{\N}{\mathbb N}

\newcommand{\R}{\mathbb R}

\newcommand{\C}{\mathbb C}

\title[Dynamics of composition operators]{Dynamics of composition operators on spaces of holomorphic functions on plane domains}
\author[J.\ B\`{e}s and C.\ Foster]{J.\ B\`{e}s and C.\ Foster}
\address{J.~B\`{e}s, Department of Mathematics and Statistics,
Bowling Green State University,
Bowling Green, Ohio 43403,
USA.}
\email{jbes@bgsu.edu}
\address{C. Foster, Department of Mathematics and Statistics,
Bowling Green State University,
Bowling Green, Ohio 43403,
USA.}
\email{fostecd@bgsu.edu}
\date{August 9th, 2024}

\subjclass{Primary 47A16; Secondary 30K20, 47B33}

\keywords{Composition operators, weighted composition operators, chaotic operators, hypercyclic operators, supercyclic operators, frequently hypercyclic operators, mixing operators}

\allowdisplaybreaks[3]
\begin{document}
\begin{abstract}We study the dynamic behaviour of (weighted) composition operators on the space of holomorphic functions on a plane domain.  Any such operator is hypercyclic if and only if it is topologically mixing, and when the symbol is automorphic, such an operator is supercyclic if and only if it is mixing. When the domain is a punctured plane, a composition operator is supercyclic if and only if it satisfies the Frequent Hypercyclity Criterion, and when the domain is conformally equivalent to a punctured disc, such an operator is hypercyclic if and only if it satisfies the Frequent Hypercyclicity Criterion. When the domain is finitely connected and either conformally equivalent to an annulus or having two or more holes, no weighted composition operator can be supercyclic. 
 \end{abstract}
\maketitle
%{\large   
\section{Introduction}
Throughout this paper, $H(\Omega)$ denotes the space of holomorphic functions on an arbitrary domain $\Omega$ of the complex plane $\mathbb{C}$, endowed with the compact open topology. Each holomorphic self map $\psi$ of $\Omega$ induces an unweighted composition operator $C_\psi: H(\Omega)\to H(\Omega)$, \ $f\mapsto f\circ \psi$, and each $\omega\in H(\Omega)$ induces an operator $M_\omega$ on $H(\Omega)$ of pointwise multiplication by the weight symbol $\omega$.
We are interested in the dynamics of the %weighted 
composition operator 
$C_{\omega, \psi}:=M_\omega C_\psi :H(\Omega)\to H(\Omega)$,\[  C_{\omega, \psi}(f)(z):=\omega (z)\ (f\circ\psi )(z) \ \ (z\in\Omega).
 \]
 We refer to $\omega$ and $\psi$ as the multiplier and symbol of $C_{\omega, \psi}$, respectively.
 %%%.  \cite{Bes2013}. \cite{BesCAOT}
 %%A characterization of hypercyclicity for the unweighted case by Grosse-Erdmann and Mortini \cite{GroMor2009} 
 Goli\'nski and Przestacki \cite{golprz2021} determined a full characterization of hypercyclicity, completing prior work by several authors on the dynamics for the unweighted case or the case when $\Omega$ is simply connected, see  \cite{GroMor2009}, \cite{Bes2013}, \cite{BesCAOT} and the references therein. A key role for this characterization is played by the holes of $\Omega$ and the dynamics of the symbol $\psi$ on the $\Omega$-convex compact subsets of $\Omega$, see Definition~\ref{D:hole}. In particular, Goli\'nski and Przestacki showed that $H(\Omega)$ supports no hypercyclic weighted composition operators whenever $\Omega$ has exactly $n$ holes with $n\ge 2$ \cite[Theorem 6.1]{golprz2021}.
 When $\Omega$ has exactly one hole, it is conformally equivalent to either the punctured plane $\mathbb{C}^*=\mathbb{C}\setminus \{ 0\}$, the punctured disc $\mathbb{D}^*=\mathbb{D}\setminus \{ 0 \}$, or an annulus $A(r)=\{ z\in \mathbb{C}: \ 1<|z|<r \}$ with $r>1$ \cite[Theorem~{10.2}]{BurckelB1979}. For the annulus case they showed that $H(A(r))$ cannot support a hypercyclic composition operator, either \cite[Theorem 6.1]{golprz2021}. For the former cases, Beltr\'an-Meneu et al  \cite[Theorem~{15}]{BeltranJordaMurillo2019} had earlier shown that no unweighted composition operator could be supercyclic on $H(\mathbb{D}^*)$ or $H(\mathbb{C}^*)$, so it was rather surprising when Goli\'nski and Przestacki found that these spaces support composition operators that are even frequently hypercyclic.
 
 \begin{proposition} {\rm ({\bf Goli\'nski-Przestacki}  \cite[Proposition~5.2]{golprz2021}) } \label{P:proposition17} 
Let $a\in \mathbb{D}^*$ and $\Omega \subseteq \mathbb{C}$ be simply connected with $0\in \Omega$ so that $a\Omega\subseteq \Omega$.  Then for each $n \in \mathbb{N}$ and  $c \neq 0$ the operator \ $C_{cz^{n}, az}: H(\Omega^{*}) \to H(\Omega^{*})$ is frequently hypercyclic.
%satisfies the Frequent Hypercyclicity Criterion, and in particular is frequently hypercyclic, mixing and chaotic.
\end{proposition}

In the same paper they posed the following.
\begin{problem} ({\bf Goli\'nski-Przestacki \cite[Problem~{5.8}]{golprz2021}})
Is every weighted composition operator acting on $H(\mathbb{C}^*)$ also frequently hypercyclic? In particular, is $C_{\omega, \frac{z}{2}}$ frequently hypercyclic on $H(\mathbb{C}^*)$, where $\omega(z)=z \exp (\frac{1}{z})$?
\end{problem}

 In this paper we continue this study and consider also the notion of supercyclicity. An earlier result by Bernal-Gonz\'alez et al~\cite{BerBonCal2007} states that for automorphic symbols of a simply connected domain of the complex plane, (unweighted) compositional supercyclicity equals hypercyclicity. In this spirit we provide the following.
\begin{theorem}\label{T:overall}
Let $\omega :\Omega\to \mathbb{C}$ and $\psi:\Omega\to \Omega$ be holomorphic, where  $\Omega$ is an arbitrary domain in $\C$. Then $C_{\omega, \psi}$ is hypercyclic on $H(\Omega)$ if and only if it is mixing. If the symbol $\psi$ is an $\Omega$-automorphism, then $C_{\omega, \psi}$ is supercyclic on $H(\Omega)$ if and only if it is mixing.
\end{theorem}

We can say more when specifying additional assumptions on $\Omega$. For the punctured plane case, we solve Problem 1 in the affirmative and show that those properties between being supercyclic to satisfying the Frequent Hypercyclicity Criterion are all equivalent, see Subsection~\ref{sub:notation} for definitions. Indeed, we have the following

\begin{theorem} \label{T:supercyclic,cstar} Let $\omega:\mathbb{C}^*\to\mathbb{C}$ and $\psi:\mathbb{C}^*\to\mathbb{C}^*$ be holomorphic, where $\mathbb{C}^*=\mathbb{C}\setminus\{ 0\}$.
The following are equivalent: 
\begin{enumerate}
\item[(1)]\ $C_{\omega, \psi}$ is supercyclic on $H(\mathbb{C}^{*})$.

\item[(2)]\ $C_{\omega, \psi}$ satisfies the Frequent Hypercyclicity Criterion on $H(\mathbb{C}^{*})$, and in particular is frequently hypercyclic, mixing and chaotic.

\item[(3)] There exist $n\in \mathbb{N}$, $W\in H(\mathbb{C}^*)$, and $a\in \mathbb{D}^*$ so that either 
\[
(i) \ \ \begin{cases} \omega(z)=z^n e^{W(z)}\\
\psi(z)=az
\end{cases}
\ \ \ \mbox{ or } \ \ \ \ \ 
(ii) \ \ \begin{cases}
 \omega(z)=z^{-n} e^{W(z)}\\
\psi(z)=\frac{1}{a}z.
\end{cases}
\]
hold.
\end{enumerate}
\end{theorem}

For the punctured disk case we have the following.

\begin{theorem} \label{T:mixing,dstar} 
Let $\omega:\mathbb{D}^*\to\mathbb{C}$ and $\psi:\mathbb{D}^*\to\mathbb{D}^*$ be holomorphic, where $\mathbb{D}^*=\mathbb{D}\setminus\{ 0\}$. The following are equivalent:
\begin{enumerate}
\item[(1)]\ $C_{\omega, \psi}$ is hypercyclic on $H(\mathbb{D}^{*})$.

\item[(2)]\ $C_{\omega, \psi}$ satisfies the Frequent Hypercyclicity Criterion on $H(\mathbb{D}^{*})$, and in particular is frequently hypercyclic, mixing and chaotic.

\item[(3)] The multiplier $\omega$ is of the form $\omega (z)=z^n e^{W(z)}$ $(z\in \mathbb{D}^*)$, where $W\in H(\mathbb{D}^*)$ and $n\in\mathbb{N}$, and the symbol $\psi $ is injective, not surjective, and satisfies $\lim_{z\to 0}\psi(z)=0$.

\end{enumerate}
\end{theorem}

Moreover, no automorphic symbol $\psi$ of $\mathbb{D^*}$ can induce a supercyclic composition operator, see Proposition~\ref{P:noscDast}. It follows that for composition operators with automorphic symbols on a punctured simply connected domain, being supercyclic is  equivalent to satisfying the Frequent Hypercyclicty Criterion, see Corollary~\ref{C:sc=FHCpunctured}.

We also note that on those cases where $\Omega$ does not allow hypercyclicity phenomena for composition operators on $H(\Omega)$, supercyclicity cannot happen for such operators, either.

\begin{theorem} \label{T:NoSChere}
Suppose $\Omega$ is finitely connected and either conformally equivalent to a (non-degenerate) annulus or has two or more holes.
Then no weighted composition operator on $H(\Omega)$ can be supercyclic.
\end{theorem}

Finally, for the case when $\Omega$ has infinitely many holes, we have that compositional supercyclicity and compositional mixing are equivalent.

\begin{theorem} \label{T:supercyclic=mixing,infinite}
Let  $\omega: \Omega \to \mathbb{C}$ and $\psi: \Omega \to \Omega$ be holomorphic, where  $\Omega \subset \mathbb{C}$ is an infinitely connected domain. The following are equivalent: 
\begin{enumerate}
\item[(1)]\ $C_{\omega, \psi}$ is supercyclic on $H(\Omega)$. 

\item[(2)]\ $C_{\omega, \psi}$ is mixing on $H(\Omega)$. 

\item[(3)] The multiplier $\omega$ is zero-free and the symbol $\psi$ is univalent, $\Omega$-convex and strongly run-away.
\end{enumerate}
\end{theorem}

The paper is organized as follows. We complete the introduction with a subsection on definitions and preliminaries. In Section~\ref{S:Nec} we establish necessary conditions for a composition operator to be supercyclic. Theorem~\ref{T:supercyclic,cstar} and Theorem~\ref{T:mixing,dstar} are proved in Section~\ref{S:C*D*}, and in Section~\ref{S:2ormore} we show Theorem~\ref{T:NoSChere}, Theorem~\ref{T:supercyclic=mixing,infinite} and Theorem~\ref{T:overall}.

\subsection{Definitions and preliminaries} \label{sub:notation}
For excellent references on linear dynamics we refer to the books by Bayart and Matheron \cite{BayMat2009} and by Grosse-Erdmann and Peris \cite{GroPer2011}. %For composition operators, we refer to the two classics by Cowen and MacCluer \cite{CowMac1995} and by Shapiro \cite{Sha1993}.
We recall basic definitions and well-known facts we use here.

\begin{definition}
An operator $T$ on a topological vector space $X$ is {\bf hypercyclic } (respectively, {\bf supercyclic}) provided for some $f\in X$ the orbit
\[
\mbox{Orb}(T, f )= \{ T^nf: \ n=0,1, \dots \}
\]
 (respectively, the set $\C\cdot \mbox{Orb}( T, f)= \{ \lambda T^nf:\ \lambda\in\C , \ n=0,1,\dots \}$)   is dense in $X$. Such $f$ is called a {\bf hypercyclic vector} for $T$ (respectively, a {\bf supercyclic vector} for $T$). When $X$ is a separable Baire space, Birkhoff's Transitivity Theorem \cite[Theorem~{1.16}]{GroPer2011}
ensures that  $T$ is hypercyclic if and only if it is topologically transitive. 
\end{definition}

\begin{definition}
An operator $T$ on a topological vector space $X$  is  {\bf topologically transitive} (respectively, {\bf mixing}) provided for each pair of non-empty open subsets $U$ and $V$ of $X$ there exists an integer $n$ so that
 \[T^n(U)\cap V\ne \emptyset\] 
 (respectively, so that $T^k(U)\cap V\ne \emptyset$ for each $k\ge n$). Also, $T$ is called {\bf weak mixing} provided the direct sum $T\oplus T$ is topologically transitive on $X\times X$, and it is said to be {\bf chaotic} provided it is topologically transitive and its set 
$
\mbox{Per}(T)=\{ f\in X: \ T^rf=f \ \mbox{ for some } r\in \N \}
$
 of periodic points is dense in $X$. 
\end{definition}
\begin{definition} %{\bf (Bayart and Grivaux)}
An operator $T$ on a separable Fr\'echet space $X$ is {\bf frequently hypercyclic} provided there is some $f\in X$ so that 
\[
\liminf_{N\to\infty} \frac{\mbox{card}\{ 0\le n\le N:\  \ T^nf\in U   \}         }{N+1}  >0
\]
for each non-empty open subset $U$ of $X$
Also, we say that $T$ satisfies the {\bf Frequent Hypercyclicity Criterion} provided there exist 
a dense subset $X_0$ of $X$ and a map $S:X_0\to X_0$ satisfying, for each $f\in X_0$, that 
\begin{enumerate}
\item[(i)]\ The series $\sum_{n=1}^\infty T^n(f)$ converges unconditionally.
\item[(ii)]\ The series $\sum_{n=1}^\infty S^n(f)$ converges unconditionally.
\item[(iii)]\ $TSf=f$.
\end{enumerate}
\end{definition}

\begin{remark} \label{R:inverseFHC}
When an operator $T$ on a separable Fr\'echet space is invertible, it is simple to see that each of the properties of being supercyclic, hypercyclic, chaotic, weak-mixing, mixing, or of satisfying the Frequent Hypercyclicity Criterion  is satisfied by $T$ if and only if it is satisfied by its inverse $T^{-1}$.  In contrast, it is possible for $T$ to be frequently hypercyclic with $T^{-1}$ not being frequently hypercyclic.
\end{remark}

%\begin{definition} \label{D:FHC}
%Let $X$ be a separable Fr\'echet space. We say that an operator $T:X\to X$ satisfies the {\bf Frequent Hypercyclicity Criterion} provided there exists 
%a dense subset $X_0$ of $X$ and a map $S:X_0\to X_0$ satisfying, for each $f\in X_0$, that 
%\begin{enumerate}
%\item[(i)]\ The series $\sum_{n=1}^\infty T^n(f)$ converges unconditionally.
%\item[(ii)]\ The series $\sum_{n=1}^\infty S^n(f)$ converges unconditionally.
%\item[(iii)]\ $TSf=f$.
%\end{enumerate}
%\end{definition}

\begin{theorem} {\rm (\cite[Theorem 2.1]{bayart_grivaux2006},\cite[Remark 2.2(b)]{bonilla_grosse-erdmann2007})}
Let $T$ be an operator on a separable Fr\'echet space and which satisfies the Frequent Hypercyclicity Criterion. Then $T$ is frequently hypercyclic, mixing and chaotic.
\end{theorem}

\begin{definition}
An operator $\tilde{T}:\tilde{X}\to \tilde{X}$ is {\bf quasi-conjugate} (respectively, {\bf linearly quasi-conjugate}) to an operator $T:X\to X$ provided there exists a continuos (respectively, continuos and linear) map $J:X\to \tilde{X}$ of dense range so that the diagram

\begin{equation}
\begin{CD}
\tilde{X}@> \tilde{T} >> \tilde{X} \\
@A J AA @ AA J A \\
X @ > T  >> X
\end{CD}
\end{equation}

\noindent
commutes. If there exists such a map $J$ that is also a homeomorphism, then we say that $\tilde{T}$ and $T$ are {\bf topologically conjugate}. Finally, if there exists such a map $J$ that is linear and a homeomorphism, then we say that $\tilde{T}$ and $T$ are {\bf linearly conjugate} or {\bf isomorphically conjugate}.
\end{definition}

\begin{remark} \label{R:1.13}
It is well known that the properties of being supercyclic, hypercyclic, weak-mixing, mixing, chaotic, or frequently hypercyclic are preserved under quasi-conjugacies. We note that the property of satisfying the Frequent Hypercyclicity Criterion is preserved under isomorphic conjugacies (and more generally, under injective linear quasi-conjugacies). Indeed, let $\tilde{T}$ and $T$ be operators on separable $F$-spaces $\tilde{X}$ and $X$ respectively and let $J:X\to \tilde{X}$ be an injective operator with dense range so that $JT=\tilde{T}J$. If $T$ satisfies the Frequent Hypercyclicity Criterion, there exists $X_0\subset X$ dense and $S:X_0\to X_0$ satisfying for each $f\in X_0$, that 
\begin{enumerate}
\item[(i)]\ The series $\sum_{n=1}^\infty T^n(f)$ converges unconditionally.
\item[(ii)]\ The series $\sum_{n=1}^\infty S^n(f)$ converges unconditionally.
\item[(iii)]\ $TSf=f$.
\end{enumerate}
Then $\tilde{X}_0:=J(X_0)$ is dense in $\tilde{X}$ and the map $\tilde{S}:\tilde{X}_0\to \tilde{X}_0$ given by $\tilde{S} (Jf):=JSf$ $(f\in X_0)$  -well defined thanks to the injectivity of $J$- satisfies for each $J(f)\in \tilde{X_0}$ and each bijection $\sigma:\mathbb{N}\to \mathbb{N}$ that
\[
\sum_{n=1}^\infty \tilde{T}^{\sigma(n)}Jf=\sum_{n=1}^\infty JT^{\sigma(n)}f=J\sum_{n=1}^\infty T^{\sigma(n)}f
\]
converges (by (i) and since $J$ is linear and continuous). Similarly, by (ii)
\[
\sum_{n=1}^\infty \tilde{S}^{\sigma(n)}Jf=\sum_{n=1}^\infty JS^{\sigma(n)}f=J\sum_{n=1}^\infty S^{\sigma(n)}f
\]
converges. But $\tilde{T}\tilde{S}Jf=\tilde{T}JSf=JTSf=Jf$ by (iii). So $\tilde{T}$ satisfies the Frequent Hypercyclicity Criterion.
\end{remark}

We conclude this introduction by recalling notions from the study of compositional dynamics on $H(\Omega)$.
\begin{definition}\label{D:hole}
Let $\Omega$ be a domain in $\mathbb{C}$.
\begin{enumerate}
\item[(i)]\ Any bounded component of  $\mathbb{C}\setminus \Omega$ is called a {\bf hole} of $\Omega$. Similarly, for any
compact subset $K$ of $\mathbb{C}$, any bounded component of $\mathbb{C}\setminus K$ is called a {\bf hole} of $K$. 
\item[(ii)]\ We say that $\Omega$  is {\bf infinitely connected} provided it has an infinite number of holes. Else it is said to be {\bf finitely-connected}; more precisely,  $\Omega$ is {\bf n-connected} whenever it has exactly $n$ holes. For the $0$-connected case (no holes), we say that $\Omega$ is {\bf simply connected}.   
\item[(iii)]\ A compact subset $K$ of  $\Omega$  is said to be {\bf $\Omega$-convex} provided each hole of $K$ contains a point of $\mathbb{C}\setminus \Omega$. Also, a holomorphic map $\psi:\Omega\to \Omega$ is {\bf $\Omega$-convex} provided $\psi(K)$ is $\Omega$-convex for each $\Omega$-convex compact subset of $\Omega$.
\item[(iv)]\ A holomorphic map $\psi:\Omega\to \Omega$ is {\bf run-away} (respectively, {\bf strongly run-away}) provided for each compact subset $K$ of $\Omega$ 
\[
\psi_n(K)\cap K=\emptyset
\]
for some $n$ (respectively, for all large $n$) in $\mathbb{N}$. Here $\psi_n$ denotes the $n$-fold composition of $\psi$ with itself.
\end{enumerate}
\end{definition}

\section{Necessary Conditions for Supercyclicity} \label{S:Nec}
We summarize in Theorem~\ref{T:necessary} below necessary conditions for supercyclicity.

\begin{theorem} \label{T:necessary}
Let $\Omega$ be a domain in $\mathbb{C}$ and let $\omega: \Omega \to \mathbb{C}$ and  $\psi: \Omega \to \Omega$ be holomorphic. If $C_{\omega, \psi}$ is supercyclic on $H(\Omega)$, the following hold:
\begin{enumerate}
\item[{\rm (i)}]\ The multiplier $\omega$ is zero-free.
\item[{\rm (ii)}]\ The symbol $\psi$ is injective, strongly run-away, and $\Omega$-convex.
\end{enumerate}
\end{theorem}

Theorem~\ref{T:necessary} follows once we establish Lemma~\ref{L:supercyclic,necessary,complex}, Lemma~\ref{L:supercyclic,convex}  and Lemma~\ref{L:supercyclic,strongrunaway} below.

\begin{lemma} \label{L:supercyclic,necessary,complex} Let $\Omega$ be a domain in $\mathbb{C}$ and $\omega: \Omega \to \mathbb{C}, \psi: \Omega \to \Omega$ be holomorphic. If $C_{\omega, \psi}$ is supercyclic on $H(\Omega)$, the following must hold:
\begin{enumerate}
\item[{\rm(i)}]\ $\omega$ is non-vanishing.

\item[{\rm (ii)}]\ $\psi$ is univalent.

\item[{\rm (iii)}]\ $\psi$ doesn't have any fixed points.

\item[{\rm (iv)}]\ For every $z \in \Omega$, the set $\overline{\{\psi_n (z): n \in \mathbb{N}_0\}}$ is not compact in $\Omega$. 
\end{enumerate}
\end{lemma}

\begin{proof}
Conclusions $(i)-(iii)$ hold under the weaker assumption that $C_{\omega, \psi}$ is weakly supercyclic \cite{BesCAOT,BesCAOTE}. To see $(iv)$, suppose that there exists $z_0\in \Omega$ so that the set
\[
K:=\overline{\{\psi_n (z_0): n \in \mathbb{N}_0\}}
\]
is compact in $\Omega$. By $(i)$, there exist $M, m>0$ so that 
\begin{equation} \label{eq:w}
m\le |\omega (z) |\le M
\end{equation}
for all $z\in K$. Since the set of supercyclic vectors of a supercyclic operator must be dense and
$U:=\{ h\in H(\Omega): \ 1< |g|< 2 \mbox{ on } K\cup \psi(K)  \}$ is open and non-empty, we may pick $f\in U$ that is supercyclic for $C_{\omega, \psi}$. Since $\{ z_0, \psi(z_0)\}$ is $\Omega$-convex, Runge's theorem ensures that there exists $g\in H(\Omega)$ satisfying
\begin{equation} \label{eq:L1}
|g(\psi(z_0))|> 2M \ \ \mbox{ and } \ \ 0 < |g(z_0)|< m.
\end{equation}
Since $f$ is supercyclic for $C_{\omega, \psi}$ there exist sequences $(n_j)$ in $\mathbb{N}_0$ and $(\lambda_j)$ in $\mathbb{C}\setminus \{ 0 \}$ so that
\[
\lambda_j C_{\omega, \psi}^{n_j}(f)\underset{j\to\infty}{\to} g \ \ \ \mbox{ in $H(\Omega)$.}
\]
In particular,
\[
\begin{aligned}
\lambda_j C_{\omega, \psi}^{n_j}(f)(\psi(z_0))&\underset{j\to\infty}{\to} g(\psi(z_0)) \\
\lambda_j C_{\omega, \psi}^{n_j}(f)(z_0)&\underset{j\to\infty}{\to} g(z_0)\ne 0,
\end{aligned}
\]
giving
\begin{equation}\label{eq:lim}
\left| \frac{ \lambda_j C_{\omega, \psi}^{n_j}(f)(\psi(z_0))}{ \lambda_j C_{\omega, \psi}^{n_j}(f)(z_0)} \right| \underset{j\to\infty}{\to} \left| \frac{ g(\psi(z_0)) }{ g(z_0)} \right| > \frac{2M}{m}.
\end{equation}
But since $f\in U$, by $\eqref{eq:w}$ and $\eqref{eq:L1}$ for each $j$ we have
\[
\begin{aligned}
\left| \frac{ \lambda_j C_{\omega, \psi}^{n_j}(f)(\psi(z_0))}{ \lambda_j C_{\omega, \psi}^{n_j}(f)(z_0)} \right| &=   \left| \frac{ f(\psi_{n_j+1}(z_0)) \prod_{i=0}^{n_j-1} \omega(\psi_i(\psi(z_0))) }{  f(\psi_{n_j}(z_0)) \prod_{i=0}^{n_j-1} \omega(\psi_i(z_0)) }  \right|   \\
&= \left| \frac{f(\psi_{n_j+1}(z_0))}{  f(\psi_{n_j}(z_0))}\right | \, \left| \frac{\omega(\psi_{n_j}(z_0))}{\omega(z_0)} \right| \\
&\le \ \ \ 2 \, \frac{M}{m}
\end{aligned}
\]
a contradiction with $\eqref{eq:lim}$.
\end{proof}

Lemma~\ref{L:supercyclic,convex} and Lemma~\ref{L:supercyclic,strongrunaway} slightly improve  \cite[Lemma~{3.3}, Lemma~{3.2}]{golprz2021} respectively where the stronger assumption of $C_{\omega, \psi}$ being hypercyclic is made. The proofs of these corresponding lemmas run parallel but use Lemma~\ref{L:supercyclic,necessary,complex}  -which is the supercyclic version of \cite[Lemma 3.1]{golprz2021} and where part (iv) is harder to establish, see \cite{BesCAOT})-.

\begin{lemma} \label{L:supercyclic,convex} Let $\Omega \subset \mathbb{C}$ be a domain which is not simply connected and let $\omega: \Omega \to \mathbb{C}$, $\psi: \Omega \to \Omega$ be holomorphic. If $C_{\omega, \psi}$ is supercyclic on $H(\Omega)$, then the map $\psi$ is $\Omega$-convex. That is, it satisfies that
for each $\Omega$-convex compact set K in $\Omega$, the set $\psi (K)$ is $\Omega$-convex.
\end{lemma}

To prove Lemma~\ref{L:supercyclic,convex}, we use the following geometrical facts.
Parts (i) and (iii) are respectively \cite[Lemma~3.10, Lemma~3.11]{GroMor2009}, and Part (ii) is sometimes referred to as the Hole Invariance Principle and holds for arbitrary compact sets~\cite[p. 276]{Remmert}.
\begin{lemma} \label{L:cons}
Let $\psi:\Omega\to \Omega$ be univalent, where $\Omega$ is a domain in $\mathbb{C}$, and let $K$ be an $\Omega$-convex compact subset of $\Omega$. Then
\begin{enumerate}
\item[(i)]\ $K$ has at most finitely many holes.
\item[(ii)]\ $\psi(K)$ has the exact number of holes as $K$ does.
\item[(iii)]\ If $L$ is a compact subset of $\Omega$ containing $K$ so that $\psi(L)$ is $\Omega$-convex, then $\psi(K)$ must be $\Omega$-convex.
\end{enumerate}
\end{lemma}

%\begin{lemma} \label{L:GE,finiteholes}{\bf \cite[Lemma 3.10]{GroMor2009}} Let K be an $\Omega$-convex compact subset of a domain $\Omega \subset \mathbb{C}$. Then K has at most finitely many holes. 
%\end{lemma}

%\begin{lemma} \label{L:GE,convexity}{\bf \cite[Lemma 3.11]{GroMor2009}} Let K and L be compact subsets of a domain $\Omega \subset \mathbb{C}$ with $K \subset L$. Let $\psi$ be a holomorphic self-map of $\Omega$ that is injective on some neighborhood of L. If K and $\psi(L)$ are $\Omega$-convex, then so is $\psi(K)$.
%\end{lemma}

\begin{proof}[Proof of Lemma~\ref{L:supercyclic,convex}] Let K be an $\Omega$-convex compact set in $\Omega$.  We want to show that $\psi(K)$ is $\Omega$-convex. By Lemma~\ref{L:supercyclic,necessary,complex} the symbol $\psi$ must be univalent, and by Lemma~\ref{L:cons}(ii) we may assume that K has at least one hole. %First we replace K by a bigger set with a better boundary.

For every $k \in \mathbb{N}$ consider the grid of points $z = x + iy$ for which x or y are integer multiples of $\frac{1}{2^k}$. Let $L_k$ be the union of closed squares with their sides lying on the grid which lie entirely in $\Omega \cap \{z \in \mathbb{C}: \lvert z \rvert < k\}$. It is clear that each $L_k$ is an $\Omega$-convex compact set. Let $k_0$ be such that the set K is contained in the interior of $L_{k_0} =: L$. If $L$ is not connected, then we may enlarge it to a connected $\Omega$-convex set (which we still denote by L) in such a way that it has a Jordan curve as its outer boundary. To prove that $\psi (K)$ is $\Omega$-convex, by Lemma~\ref{L:cons}(iii) it suffices to show that $\psi (L)$ is $\Omega$-convex.

By Lemma~\ref{L:cons}(i)  the compact set $L$ has finitely many holes, and by construction each hole of $K$ contains a hole of $L$. So $L$ has $p$ holes,  for some $p\ge 1$. Also, the boundary of $L$ consists of
(piece-wise smooth) Jordan curves $\gamma_0, ..., \gamma_p$, where $\gamma_0$ forms the outer boundary of L. We assume that $\gamma_0$ is negatively oriented and $\gamma_1, ..., \gamma_p$ are positively oriented. That is, for each $\ell=0,\dots, p$ the set $L$ is to the left of $\gamma_\ell$ as we traverse $\gamma_\ell$ with its assigned orientation.

By means of contradiction, assume that  $\psi (L)$ is not $\Omega$-convex. So there exists a hole $\mathcal{O}$ of $\psi (L)$ which does not contain a point of $\mathbb{C} \setminus \Omega$. Since univalent maps send boundaries to boundaries, there exists $0 \leq j \leq p$ such that the curve $\psi \circ \gamma_j$ is the boundary of $\mathcal{O}$. Observe that $\overline{\mathcal{O}}\subset \Omega$ and, since univalent maps preserve orientation, $\psi(L)$ lies on the left of  
$\psi \circ \gamma_j$ and thus  $\mathcal{O}$ lies on the right of  
$\psi \circ \gamma_j$ as we traverse the oriented Jordan curve $\psi \circ \gamma_j$. We show that this leads to a contradiction. 

Let $z_0$ be a point from $\mathbb{C} \backslash \Omega$ which is contained in the hole of $\gamma_j$. Let $g\in H(\Omega)$, $g(z)=(z-z_0)^\alpha$, where  $\alpha$ is the integer
%\[
\begin{equation} \label{eq:alfa}
 \alpha= \begin{cases}
  \ \ 1+\frac{1}{2\pi i}\int_{\gamma_j} \frac{\omega'(z)}{\omega(z)} dz & \mbox{if } j = 0 \\
  -1- \frac{1}{2\pi i}\int_{\gamma_j} \frac{\omega'(z)}{\omega(z)} dz & \mbox{if } j \in\{ 1, ..., p\}.
\end{cases}
\end{equation}
%\]
if $\omega$ is non-constant. If $\omega$ is constant,  define $\alpha:=1$ if $j=0$ and $\alpha:=-1$ if ${if } j \in\{ 1, ..., p\}.$
By the orientations chosen for the $\gamma_\ell$'s the Argument Principle gives
\[
\frac{1}{2 \pi i} \int_{\gamma_j} \frac{g'(z)}{g(z)} dz = 
\begin{cases}
 \ \   \alpha & \mbox{if } j = 0 \\
  -\alpha & \mbox{if } j \in\{ 1, ..., p\}.
\end{cases}
 \]
Since $C_{\omega, \psi}$ is supercyclic, there exists $f \in H(\Omega)$ and sequences $(n_k)$ in $\mathbb{N}_0$ and $(\lambda_k)$ in $\mathbb{C}$ such that 
\[
\lambda_k C_{\omega, \psi}^{n_k} (f)\underset{k\to\infty}{\to} g \ \ \ \mbox{ in $H(\Omega)$.}
\]
Moreover, we may assume $f\in H(\Omega)\setminus \mbox{span}\{ g \}$ and thus that $n_k\ge 1$ for all $k$.
Since $g$ is zero-free on a neighborhood of $L$, we may assume by 
Hurwitz' theorem that $\lambda_k C_{\omega, \psi}^{n_k} (f)$ is zero-free on $L$ (and in particular $f$ is zero-free on $\psi_{n_k}\circ \gamma_j$ ) for large $k$ and by the continuity of the operator of complex differentiation on $H(\Omega)$ that
\begin{equation}
\lim_{k\to\infty} \frac{1}{2 \pi i} \int_{\gamma_j} \frac{(\lambda_{k} C_{\omega, \psi}^{n_k} (f))' (z)}{\lambda_k C_{\omega, \psi}^{n_k} (f)(z)} dz = \frac{1}{2 \pi i} \int_{\gamma_j} \frac{g'(z)}{g(z)} dz.
\end{equation}

Now, notice that for each $m\in \mathbb{N}$ the closure of the hole $\psi_{m-1}(\mathcal{O})$ of $\psi_m\circ \gamma_j$ lies  in $\Omega$ and by Cauchy's theorem
\begin{equation} \label{eq:a}
\int_{\psi_m \circ \gamma_j} \frac{\omega' (z)}{\omega (z)} dz=0
\end{equation}
since $\omega$ is holomorphic and zero-free on $\Omega$. 
Moreover, for large $k\in \mathbb{N}$ the hole of $\psi_{n_k}\circ \gamma_j$ is contained in $\Omega$ and 
to the right of $\psi_{n_k}\circ \gamma_j$ and since
$f$ is holomorphic on a neighborhood of the closure of this hole and is zero free on $\psi_{n_k}\circ \gamma_j$, by the Argument Principle and the orientation of $\psi_{n_k}\circ \gamma_j$ we must have
\begin{equation} \label{eq:b}
\frac{1}{2\pi i} \int_{\psi_{n_k} \circ \gamma_j} \frac{f' (z)}{f(z)} dz \le 0.
\end{equation}
Thus by $\eqref{eq:a}$  for large $k$ we have
\[
\begin{aligned}
\int_{\gamma_j} \frac{(\lambda_{k} C_{\omega, \psi}^{n_k} (f))^{'} (z)}{\lambda_{k}C_{\omega, \psi}^{n_k} (f)(z)} dz & = \sum_{m=0}^{n_k-1} \int_{\gamma_j} \frac{\omega^{'} (\psi_m (z)) \psi_m^{'} (z)}{\omega(\psi_m (z))} dz + \int_{\gamma_j} \frac{f^{'} (\psi_{n_k} (z)) \psi_{n_k}'(z)}{f(\psi_{n_k}(z))} dz &\\
& = \sum_{m=0}^{n_k-1} \int_{\psi_m\circ \gamma_j} \frac{\omega' (z)}{\omega (z)} dz + \int_{\psi_{n_k} \circ \gamma_j} \frac{f' (z)}{f(z)} dz \\
&= \int_{\gamma_j} \frac{\omega' (z)}{\omega (z)} dz + \int_{\psi_{n_k} \circ \gamma_j} \frac{f' (z)}{f(z)} dz. 
\end{aligned}
\]

\noindent It follows by $\eqref{eq:alfa}$ that 
\[
%\begin{aligned}
\lim_{k\to\infty} \frac{1}{2 \pi i} \int_{\psi_{n_k} \circ \gamma_j} \frac{f' z)}{f(z)} dz = \frac{1}{2 \pi i} \int_{\gamma_j} \frac{g'(z)}{g(z)} dz -  \frac{1}{2 \pi i}  \int_{\gamma_j} \frac{\omega' (z)}{\omega (z)} dz 
= 1,
%\end{aligned}
\]
a contradiction with $\eqref{eq:b}$
\end{proof}

\begin{lemma} \label{L:supercyclic,strongrunaway} Let $\Omega \subset \mathbb{C}$ be a domain and let $\omega: \Omega \to \mathbb{C}$ and $\psi: \Omega \to \Omega$ be holomorphic. If $C_{\omega, \psi}$ is supercyclic on $H(\Omega)$, then $\psi$ has the strong run-away property.

\end{lemma}

\begin{proof} We have 3 cases.

\noindent Case 1: $\Omega = \mathbb{C}$. By Lemma \ref{L:supercyclic,necessary,complex} the symbol $\psi: \mathbb{C} \to \mathbb{C}$ is univalent and without fixed points, so there exists $0 \neq b \in \mathbb{C}$ such that $\psi (z) = z+b$ and $\psi$ is strongly run-away.

\noindent Case 2: $\mathbb{C} \backslash \Omega$ is a singleton. We may assume $\Omega = \mathbb{C}^*=\mathbb{C} \backslash \{0\}$. Since $\psi: \mathbb{C}^*\to \mathbb{C}^*$ is univalent, there exists $0 \neq a \in \mathbb{C}$ such that $\psi (z) = az$ or $\psi(z) = \frac{a}{z}$. But the latter has a fixed point, and by Lemma~\ref{L:supercyclic,necessary,complex}(iv) we must have $\lvert a \rvert \ne 1$. That is, the supercyclicity of $C_{\omega, \psi}$ implies that $\psi (z) = az$ with $\lvert a \rvert \neq 0, 1$ forcing $\psi$ to be strongly run-away.

\noindent Case 3: $\mathbb{C} \backslash \Omega$ contains at least two points. Here $\Omega$ is a hyperbolic Riemann surface \cite[Lemma~2.5]{Milnor2006}. 
Lemma \ref{L:supercyclic,necessary,complex} and \cite[Theorem 5.2]{Milnor2006} now give that $\psi$ is strongly run-away.
\end{proof}

\section{The punctured plane and punctured disc cases} \label{S:C*D*}
\subsection{Preliminaries on punctured simply connected domains}
Throughout this subsection, $\Omega \subseteq \mathbb{C}$ is a simply connected domain containing zero and $a\in\mathbb{D}^*$ so that $a\Omega\subseteq \Omega$. The main goal is to establish Proposition~\ref{P:proposition19,notsupercyclic}, which is used in the next two subsections to
show Theorem~\ref{T:supercyclic,cstar}  and Theorem~\ref{T:mixing,dstar} .

We recall the following characterization of zero-free analytic maps on punctured simply connected domains used in \cite{golprz2021}, see \cite{Radstrom1953} and \cite[p.15]{Bhattacharyva69}.

\begin{lemma} \label{L:5.5} \cite{Radstrom1953}
Let $\Omega$ be a simply connected domain containing zero and $\omega\in H(\Omega^*)$ be zero-free. Then there exists a unique integer $k$ and $W\in H(\Omega^*)$ so that
\[
\omega(z)=z^k \, e^{W(z)} \ \ \ \ \ \ (z\in\Omega^*).
\]
Indeed, $k$ is the winding number around zero of the curve $\omega \circ \gamma$ for all positively oriented circles $\gamma$ in $\Omega$ around zero, that is $k=\frac{1}{2\pi i} \int_{\omega \circ \gamma} \frac{1}{\eta} d\eta$.
\end{lemma}

\begin{lemma} \label{L:corollary15}
\textup{\cite[Proposition~5.1]{golprz2021}}. 
Let $\Omega \subseteq \mathbb{C}$ be simply connected containing zero and $a\in\mathbb{D}^*$ so that $a\Omega\subseteq \Omega$,
and let $\omega \in H(\Omega)$ be zero-free. Then 
\begin{enumerate}
\item[(i)] The infinite product $h:\Omega\to \C$, $h(z):=\prod_{j=0}^\infty \frac{\omega(a^jz)}{\omega(0)}$ converges locally uniformly in $\Omega$. Moreover, $h\in H(\Omega)$ is a zero-free eigenfunction for $C_{\omega, az}$ of eigenvalue $\omega(0)$.
\item[(ii)]  For every $\tilde{\omega} \in H(\Omega^*)$ the diagram
\[
\begin{CD}
H(\Omega^*)@>  C_{\omega(0)\tilde{\omega}, az}   >>H(\Omega^*) \\
@V M_h VV @ VV M_h V \\
H(\Omega^*) @ >  C_{\omega\tilde{\omega}, az}  >> H(\Omega^*).
\end{CD}
\]
commutes. Since $h$ is zero-free, $C_{\omega \tilde{\omega}, az}: H(\Omega^{*}) \to H(\Omega^{*})$ is isomorphically conjugate to the  $C_{\omega(0) \tilde{\omega}, az}: H(\Omega^{*}) \to H(\Omega^{*})$.

\end{enumerate}
  
\end{lemma}

\begin{lemma} \label{L:restriction} 
Let  $a\in\mathbb{D}^*$ and $\Omega \subseteq \mathbb{C}$ be simply connected and containing zero so that $a\Omega\subseteq \Omega$.
Then for any $\omega\in H(\mathbb{C}^*)$ the operator $C_{\omega, az}: H(\Omega^*) \to H(\Omega^*)$ is linearly quasi-conjugate to   $C_{\omega, az}: H(\mathbb{C^{*}}) \to H(\mathbb{C^{*}})$.
\end{lemma}

\begin{proof} For any $f \in H(\mathbb{C^{*}})$, let $r(f)$ denote the restriction of f to $\Omega^{*}$. Then the diagram
\[
\begin{CD}
H(\mathbb{C}^*)@> C_{\omega, az} >>H(\mathbb{C}^*) \\
@V r VV @ VV r V \\
H(\Omega^*) @ > C_{\omega, az} >> H(\Omega^*).
\end{CD}
\]
commutes and the restriction map $r:H(\mathbb{C}^*)\to H(\Omega^*)$ is linear and continuous.% and has dense range. To see the latter,
But $r$ has dense range too: given $g \in H(\Omega^{*})$ and $K\subset \Omega^*$ compact, by Runge's theorem, there exists a rational function $h$ with a pole at 0 that is arbitrarily close to $g$ on $K$, so $r(h)\in H(\Omega^*)$ is arbitrarily close to $g$ on $K$. 
\end{proof}

\begin{lemma} \label{L:constant,notsupercyclic} 
Let $\Omega$ be a simply connected domain containing zero and $a\in \mathbb{D}^*$ so that $a\Omega\subseteq \Omega$.
For every $b \in \mathbb{C}$, the operator $C_{b, az}$ is not supercyclic on $H(\Omega^{*})$.
\end{lemma}

\begin{proof} %Theorem~\ref{T:Jorda} tells us that the operator $C_{az}$ is not supercyclic on $H(\mathbb{C}^{*})$.
Given $f \in H(\mathbb{C}^{*})$, then the projective orbit of f under $C_{b, az}$ is contained in the projective orbit of f under $C_{az}$ since b is a constant. The latter is not dense in $H(\mathbb{C}^{*})$ by \cite[Theorem~{15}]{BeltranJordaMurillo2019}, so $C_{b, az}$ is not supercyclic on $H(\mathbb{C}^{*})$. 
\end{proof}

%\begin{lemma} \label{L:omeganeq0} 
%Let $\Omega$ be a simply connected domain containing zero and $a\in \mathbb{D}^*$ so that $a\Omega\subseteq \Omega$.
%Let $\omega \in H(\Omega)$ be such that $\omega(0) \neq 0$. Then $C_{\omega, az}$ is not supercyclic on $H(\Omega^{*})$.
%\end{lemma}

%\begin{proof} If $\omega$ has a zero in $\Omega^*$, then $C_{\omega, az}$ is not supercyclic by Lemma \ref{L:supercyclic,necessary,complex}. So assume that $\omega$ has no zeros in $\Omega$. Then by Lemma \ref{L:corollary15}, the operator $C_{\omega, az}$ is conjugate to the operator $C_{\omega(0), az}$ which is not supercyclic by Lemma \ref{L:constant,notsupercyclic}.
%\end{proof}

\begin{proposition} \label{P:proposition19,notsupercyclic}
Let $\Omega$ be a simply connected domain containing zero and $a\in \mathbb{D}^*$ so that $a\Omega\subseteq \Omega$.
 Let $k_0 \in \mathbb{N}$ and $c \neq 0$. Then $C_{{\frac{c}{z^{k_0}}}, az}$ is not supercyclic on $H(\Omega^{*})$. 
\end{proposition}

We will show Proposition~\ref{P:proposition19,notsupercyclic} by creating a quasi-conjugacy with a suitable pseudo-shift on an $F$-sequence space and using a characterization of supercyclicity for weighted pseudo-shifts due to Liang and Zhou, see Theorem~\ref{L:sequencespace} below.
%Indeed, since $0\in \Omega$ there exists $0<r<1$ so that $r\mathbb{D}^*\subset \Omega^*$, so each $f\in H(\Omega^*)$ has a unique Laurent expansion about 0 that converges on $r\mathbb{D}^*$.  If we look at the coefficients in this expansion, we can think of f as a bilateral sequence and the action of the operator $C_{\frac{c}{z^{k_0}}}, az$ as a weighted pseudo-shift. 

First we recall the notion of an $F$-sequence space and of a weighted pseudo-shift.  For a countable index set $I$, a sequence space over $I$ is any linear subspace of the space $w(I)=\mathbb{K}^I$ of all scalar sequences $(x_i)_{i\in I}$. The space $w(I)$ is endowed with the product topology. By $e_i$ $(i\in I)$ we denote the canonical unit vectors $e_i=(\delta_{i,j})_{j\in I}$. A topological sequence space over $I$ is a sequence space $X$ over $I$ that is endowed with a linear topology for which the inclusion mapping $i:X\to w(I)$ is continuous (equivalently, so that each coordinate functional $f_s:X\to \mathbb{K}$, $(x_i)_{i\in I}\mapsto x_s$, is continuous). When $X$ becomes an $F$-space (respectively, Fr\'echet space) under this topology, it is called an $F$-sequence space (respectively, Fr\'echet sequence space).

The family $( e_i )_{i\in I}$ of unit vectors is called an OP-basis (or Ovsepian-Pelczy\'nski basis) provided its linear span is dense in $X$ and its family $(f_i)_{i\in I}$ of coordinate functionals is equicontinuous. If $(e_i)_{i\in I}$ is a Schauder basis of $X$ (i.e. if every $x\in X$ has a unique representation as $x=\sum_{i\in I} x_i e_i$ with $x_i\in \mathbb{K}$ for $i\in I$ ) and $X$ is an $F$-space, then $(e_i)_{i\in I}$ is an OP-basis~\cite[Corollary 2.6.2]{RolewiczB1984}.
         
\begin{definition} Let X and Y be topological sequence spaces over indexing sets I and J, respectively. Then a continuous linear operator T: $X \to Y$ is called a weighted pseudo-shift if there is a sequence $(b_j)_{j \in J}$ of nonzero scalars and an injective mapping $\phi: J \to I$ such that
\[ T(x_i)_{i \in I} = (b_j x_{\phi(j)})_{j \in J}\]
 for $(x_i) \in X$. We write $T = T_{b, \phi}$, and $(b_j)_{j \in J}$ is called the weight sequence. 
\end{definition}

We'll also need the inverse of $\phi$ which we'll call $\psi: \phi(J) \to J$. Note that we set $b_{\psi(i)} = 0$ and $e_{\psi(i)} = 0$ if $i \in I \backslash \phi(J)$. We make use of the following characterization of supercyclicity.

\begin{theorem} \label{L:sequencespace} {\rm {\bf Liang and Zhou} \cite[Theorem 5.5]{Liang2014})}  Let X be an F-sequence space over I in which $(e_i)_{i \in I}$ is an OP-basis. Let $T = T_{b, \phi}: X \to X$ be a weighted pseudo-shift. Then the following are equivalent:
\begin{enumerate}
\item[{\rm (a)}]\ T has a dense set of supercyclic vectors.

\item[{\rm (b)}]
\begin{enumerate}
	\item[{\rm (i)}]\ The mapping $\phi: I \to I$ has no periodic points.

\item[{\rm (ii)}]\ There exists an increasing sequence $(n_k)$ of positive integers such that for every  $i, j \in I$,
\[
\left \| (\prod_{v=0}^{n_k - 1} b_{\phi^{v} (j)})^{-1} e_{\phi^{n_k} (j)}  \right \|  \,  \left \| (\prod_{v=1}^{n_k} b_{\psi^{v} (i)}) e_{\psi^{n_k} (i)} \right \|  \underset{k\to\infty}{\to} 0,
\]
where $\| \, \|$ denotes the $F$-norm of $X$.
\end{enumerate}
\end{enumerate}
\end{theorem}

Now, notice that for any $r>0$ the set $X$ of all bilateral sequences $(a_k)_{k\in \mathbb{Z}}$ for which the series $\sum_{k\in\mathbb{Z}} a_k z^k$ converges on $r\mathbb{D}^*$ is a sequence space over $\mathbb{Z}$ for which
\[
\Phi: H(r\mathbb{D}^*)\to X, \ f=\sum_{k\in \mathbb{Z}} a_k z^k \mapsto (a_k)_{k\in \mathbb{Z}} 
\]
is a linear bijection. Since $H(r\mathbb{D}^*)$ is a Fr\'echet space, we may endow $X$ with the topology induced by $\Phi$, so that $\Phi$ becomes a Fr\'echet isomorphism. Notice that for each $s\in\mathbb{Z}$, the map
\[
\delta_s:H(r\mathbb{D}^*)\to \mathbb{C},\ f=\sum_{k\in\mathbb{Z}}a_kz^k\mapsto a_s=\frac{1}{2\pi i} \int_{C(0, \frac{r}{2}} \frac{f(z)}{z^{s+1}} dz
\]
is continuous and corresponds via $\Phi$ to the s-coordinate functional on $X$, so $X$ is a Fr\'echet sequence space.  Notice also that $(z^k)_{k\in\mathbb{Z}}$ is a Schauder basis for $H(r\mathbb{D}^*)$ and that $\Phi(z^k)=e_k (k\in\mathbb{Z})$, so $(e_k)_{k\in \mathbb{Z}}$ is a Schauder basis for $X$ and thus it is an OP-basis for $X$.

We use the following notation for $H(r\mathbb{D}^*)$.
For each $m\in \mathbb{N}$ let
\[K_m := \{z \in r\mathbb{D}^{*}: \frac{r}{m} \leq \lvert z \rvert \leq r (1 - \frac{1}{m})\}\]
and
\[
\rho_m (f) = \textrm{sup}\{\lvert f(z) \rvert: z \in K_m\}. 
\]

\noindent Then $(K_m)_{m \in \mathbb{N}}$ is an exhaustion of compact sets for $r\mathbb{D}^{*}$ and
 $(\rho_m)_{m \in \mathbb{N}}$ is an increasing and separating sequence of seminorms on $H(r\mathbb{D}^{*})$ and
 \[
 \| f \|:= \sum_{m=1}^\infty \frac{1}{2^m} \mbox{min}\{ 1, \rho_m(f)\}
 \]
 defines an $F$-norm on $H(r\mathbb{D}^*)$.

\begin{proof}[Proof of Proposition~\ref{P:proposition19,notsupercyclic}]\ 
Let $0<r<1$ so that $r\mathbb{D}\subset \Omega$. By Lemma~\ref{L:restriction}, it suffices to show that the operator $T := C_{\frac{c}{z^{k_0}}, az}:H(r\mathbb{D}^{*})\to H(r\mathbb{D}^{*})$, which is quasi-conjugate to  $C_{\frac{c}{z^{k_0}}, az}:H(\Omega^{*})\to H(\Omega^{*})$, is not supercyclic. As mentioned previously, any $f \in H(r\mathbb{D}^{*})$ has a Laurent series expansion about 0, say $f = \sum_{-\infty}^{\infty} d_j z^{j}$. The operator T acts on f in the following way:

\begin{center} $(Tf)(z) = \sum_{-\infty}^{\infty} (c d_j a^{j}) z^{j-k_0} = \sum_{-\infty}^{\infty} (c d_{j+k_0} a^{j+k_0}) z^{j}$. 
\end{center}

\noindent Thus, we see that T is a weighted pseudo-shift with $b = (ca^{j+k_0})_{j \in \mathbb{Z}}$ and $\phi(j) = j + k_0$, and as discussed earlier we can view $H(r\mathbb{D}^{*})$ as an F-sequence space over $\mathbb{Z}$ and where the functions $(z^j)_{j \in \mathbb{Z}}$ form an OP-basis for $H(r\mathbb{D}^*)$. 
%We will show that condition (ii) in Lemma \ref{L:sequencespace} does not hold, so T doesn't have a dense set of supercyclic vectors. Because every %supercyclic operator has a dense set of supercyclic vectors, this will show that T is not supercyclic.  

Suppose that $T$ is supercyclic. By Theorem~\ref{L:sequencespace}, there exists a exists an increasing sequence $(n_k)$ of positive integers such that, for every $i, j \in I$ we have
\begin{equation} \label{eq:pseudoshift}
\begin{aligned}& \ \ \ \ \ \   \lvert \lvert (\prod_{v=0}^{n_k - 1} b_{\phi^{v} (j)})^{-1} e_{\phi^{n_k} (j)}  \rvert \rvert \
\lvert \lvert (\prod_{v=1}^{n_k} b_{\psi^{v} (i)}) e_{\psi^{n_k} (i)} \rvert \rvert  
=  \\
&= \lvert \lvert (\prod_{v=0}^{n_k - 1} c a^{j+k_0(v+1)})^{-1} z^{j+n_k k_0} \rvert \rvert \
\lvert \lvert (\prod_{v=1}^{n_k} c a^{i+k_0 (1-v)}) z^{i - n_k k_0} \rvert \rvert \\ 
&= \lvert \lvert p_k ^ {-1} z^{j+n_k k_0} \rvert \rvert  \
\lvert \lvert q_k z^{i-n_k k_0} \rvert \rvert \underset{k\to\infty}{\to} 0,
\end{aligned}
\end{equation}

\noindent where $p_k = \prod_{v=0}^{n_k - 1} c a^{j+k_0(v+1)}$ and $q_k = \prod_{v=1}^{n_k} c a^{i+k_0 (1-v)}$ for each $k\in \mathbb{N}$.  Note that $j + n_k k_0 \geq 1$ for k large regardless of what j is. Thus, for k large we have
\[
\rho_m (z^{j+n_k k_0}) = \textrm{sup}_{z \in K_m} \lvert z \rvert ^{j+n_k k_0} = (r(1-\frac{1}{m}))^{j+n_k k_0},
\]
and similarly
\[
\rho_m (z^{i-n_k k_0}) = (\frac{r}{m})^{i-n_k k_0}.
\]
\noindent So for k large we have
\[
\begin{aligned}
\lvert \lvert p_k ^ {-1} z^{j+n_k k_0} \rvert \rvert &= \sum_{m=1}^{\infty} \frac{1}{2^{m}} \textrm{min}\{ 1, \lvert p_k ^{-1} \rvert  (r(1- \frac{1}{m}))^{j+n_k k_0}\} \\
&\geq \frac{1}{4} \textrm{min}\{ 1, \lvert p_k ^{-1} \rvert (\frac{r}{2})^{j+n_k k_0}\}.
\end{aligned}
\]

Looking back at the definition of $p_k$, we see that most of the exponents inside the product will be positive for large k. Because $0 < \lvert a \rvert < 1$, $\lvert p_k ^{-1} \rvert$ tends to infinity as k tends to infinity. Although $(\frac{r}{2})^{j+n_k k_0}$ tends to zero, we have that $\lvert p_k ^{-1} \rvert (\frac{r}{2})^{j+n_k k_0}$ tends to infinity because the highest exponent in $p_k$ is $j+n_k k_0$, multiplied by many more positive exponents as k gets large. Thus, 
\[
\frac{1}{4} \textrm{min}\{ 1, \lvert p_k ^{-1} \rvert (\frac{r}{2})^{j+n_k k_0}\} \to \frac{1}{4} \; \textrm{as} \; k \to \infty.
\]
\noindent Similarly, for k large,
\[
\lvert \lvert q_k z^{i-n_k k_0} \rvert \rvert = \sum_{m=1}^{\infty} \frac{1}{2^{m}} \textrm{min}\{ 1, \lvert q_k \rvert (\frac{r}{m})^{i-n_k k_0}\} \geq \frac{1}{4} \textrm{min}\{1, \lvert q_k \rvert (\frac{r}{2})^{i-n_k k_0}\}.
\]

\noindent Both $\lvert q_k \rvert$ and $(\frac{r}{2})^{i-n_k k_0}$ tend to infinity as k grows large. Thus, their product tends to infinity, so we have
\[
\frac{1}{4} \textrm{min}\{ 1, \lvert q_k \rvert (\frac{r}{2})^{i-n_k k_0}\}  \to \frac{1}{4} \; \textrm{as} \; k \to \infty.
\]
\noindent Thus,
\[
\lvert \lvert p_k ^ {-1} z^{j+n_k k_0} \rvert \rvert \lvert \lvert q_k z^{i-n_k k_0} \rvert \rvert \nrightarrow 0 \; \textrm{as} \; k \to \infty,
\]
a contradiction with $\eqref{eq:pseudoshift}$. So  T isn't supercyclic on $H(r\mathbb{D}^{*})$ nor on $H(\Omega^{*})$.
\end{proof}

%
%
%.    SUBSECTION
%
%
%

\subsection{The punctured plane case. Proof of Theorem~\ref{T:supercyclic,cstar}}%. {T:supercyclic,cstar} 
\setcounter{section}{1} 
\setcounter{theorem}{2}
%\numberwithin{equation}{section} 

\begin{theorem}  %\label{T:supercyclic,cstar} 
%{\bf Theorem~\ref{T:supercyclic,cstar}.} 
Let $\omega:\mathbb{C}^*\to\mathbb{C}$ and $\psi:\mathbb{C}^*\to\mathbb{C}^*$ be holomorphic, where $\mathbb{C}^*=\mathbb{C}\setminus\{ 0\}$.
Then the following are equivalent: 
\begin{enumerate}
\item[(1)]\ $C_{\omega, \psi}$ is supercyclic on $H(\mathbb{C}^{*})$.

\item[(2)]\ $C_{\omega, \psi}$ satisfies the Frequent Hypercyclicity Criterion on $H(\mathbb{C}^*)$.

\item[(3)]\ There exist  $a\in \mathbb{D}^*$, $n\in \mathbb{N}$  and $W\in H(\mathbb{C}^*)$ so that either 
\[
(i) \ \ \begin{cases} \omega(z)=z^n e^{W(z)}\\
\psi(z)=az
\end{cases}
\ \ \ \mbox{ or } \ \ \ \ \ 
(ii) \ \  \begin{cases}
 \omega(z)=z^{-n} e^{W(z)}\\
\psi(z)=\frac{1}{a}z
\end{cases}
\]
hold.
\end{enumerate}
\end{theorem}

\begin{proof}[Proof of Theorem~\ref{T:supercyclic,cstar}] 
The implication $(2)\Rightarrow (1)$ is immediate.  $(1)\Rightarrow (3)$. Suppose  $C_{\omega, \psi}$ is supercyclic on $H(\mathbb{C}^*)$. By Theorem~\ref{T:necessary} the multiplier $\omega$ is zero-free and the symbol $\psi:\mathbb{C}^*\to \mathbb{C}^*$ is univalent and run-away. This forces $\psi$ to be of the form 
\[
\psi(z)=az \ \ (z\in \mathbb{C}^*)
\]
for some $a\in\mathbb{C}$ with $|a|\ne 0,1$. Also, since $\omega:\mathbb{C}^*\to \mathbb{C}$ is zero-free by Lemma~\ref{L:5.5} it must be of the form
\[
\omega(z)=z^k e^{f(z)+g(\frac{1}{z})} \ \ (z\in \mathbb{C}^*)
\]
for some $k\in \mathbb{Z}$ and $f, g\in H(\mathbb{C})$ with $f(0)=0$. 
Suppose first that $0<|a|<1$.  We show that Condition (3)(i) must hold and that in this case $C_{\omega, \psi}$ satisfies the Frequent Hypercyclicity Criterion. 
%Let $f, g\in H(\mathbb{C})$ with $f(0)=0$ so that
%\[
%W(z)=f(z)+g(\frac{1}{z}) \ \ \ (z\in\mathbb{C}^*).
%\]
By Lemma~\ref{L:corollary15}, the operator $C_{\omega, \psi}=C_{z^k e^{f(z)+g(\frac{1}{z})}, az}$ is isomorphically conjugate to 
\[
T_1=C_{z^k e^{g(\frac{1}{z})}, az},
\]
which in turn is isomorphically conjugate to 
\[
\begin{aligned}
T_2:=C_\frac{1}{z} T_1 C_{\frac{1}{z}}^{-1}&= C_\frac{1}{z} M_{z^k e^g(\frac{1}{z})} C_{az} C_\frac{1}{z} \\
&= M_{(\frac{1}{z})^k e^{g(z)}} C_\frac{1}{z} C_{az} C_\frac{1}{z}=  M_{(\frac{1}{z})^k e^{g(z)}} C_\frac{z}{a}.
\end{aligned}
\]
We want to show that $k\in\mathbb{N}$ and that in this case $T_2$ satifies the Frequent Hypercyclicity Criterion.
Notice first that $T_2$ is supercyclic if and only if its inverse $T_2^{-1}$ is supercyclic. 
But by Lemma~\ref{L:corollary15}
\[
\begin{aligned}
T_2^{-1}&=\left(  M_{(\frac{1}{z})^k e^{g(z)}} C_\frac{z}{a} \right)^{-1} \\
&=C_{az} M_{z^k e^{-g(z)}} =M_{(az)^k e^{-g(az)}} C_{az}
\end{aligned}
\]
is isomorphically conjugate on $H(\mathbb{C}^*)$ to
\[
T_3:=M_{(az)^k e^{-g(0)}} C_{az}=C_{cz^k, az},
\]
where $c=a^ke^{-g(0)}\ne 0$. So $T_3$ is supercyclic. By %Theorem~\ref{T:Jorda}
Lemma~\ref{L:constant,notsupercyclic}  we must have $k\ne 0$, and by Proposition~\ref{P:proposition19,notsupercyclic} we cannot have $k<0$. So $k\in \mathbb{N}$, and Condition (3)(i) holds.
Consider now the case $|a|>1$.  
\[
C_{\omega, \psi}=M_{z^ke^{f(z)+g(\frac{1}{z})}} C_{az}
\]
is isomorphically conjugate to 
\[
%\begin{aligned}
C_{\frac{1}{z}} M_{z^ke^{f(z)+g(\frac{1}{z})}} C_{az} C_{\frac{1}{z}}^{-1}=M_{z^{-k} e^{f(\frac{1}{z}) + g(z)}} C_{\frac{1}{z}}  C_{az} C_{\frac{1}{z}} =M_{z^{-k} e^{f(\frac{1}{z}) + g(z)}} C_{\frac{z}{a}},
%\end{aligned}
\]
which by Lemma~\ref{L:corollary15} is isomorphically conjugate to
\[
T_4:= M_{z^{-k}e^{f(\frac{1}{z})+g(0)}} C_{\frac{z}{a}}.
\]
So $T_4$ is supercyclic, and hence its inverse $T_4^{-1}$ is. But 
\[
T_4^{-1}=C_{az}M_{z^k e^{-f(\frac{1}{z})-g(0)}}=M_{(az)^k e^{-f(\frac{1}{az})-g(0)}} C_{az}
\]
is isomorphically conjugate  to
\[
\begin{aligned}
C_{\frac{1}{z}} T_4^{-1} C_{\frac{1}{z}}^{-1}&=C_{\frac{1}{z}}M_{(az)^k e^{-f(\frac{1}{az})-g(0)}} C_{az}C_{\frac{1}{z}} \\
&=M_{a^kz^{-k} e^{-f(\frac{z}{a})-g(0)}} C_\frac{1}{z}C_{az}C_{\frac{1}{z}}= M_{a^kz^{-k} e^{-f(\frac{z}{a})-g(0)}} C_{\frac{z}{a}},
\end{aligned}
\]
which by Lemma~\ref{L:corollary15} is isomorphically conjugate to 
\[
T_5:= M_{\tilde{c} z^{-k}} C_{\frac{z}{a}} = C_{\tilde{c} z^{-k}, \frac{1}{a} z},
\]
where $\tilde{c}=a^k e^{-f(0)-g(0)}\ne 0$.
Since $T_5$ is supercyclic, by %Theorem~\ref{T:Jorda} 
Lemma~\ref{L:constant,notsupercyclic} and Proposition~\ref{P:proposition19,notsupercyclic} we must now have $k<0$. So $(3)(ii)$ holds.

$(3)\Rightarrow (2)$ Let $f, g\in H(\mathbb{C})$ with $f(0)=0$ so that
\[
W(z)=f(z)+g(\frac{1}{z}) \ \ \ (z\in\mathbb{C}^*).
\]
Suppose first that Condition (3)(i) holds. Performing the same isomorphic conjugacies as in $(1)\Rightarrow (3)$ for the case $0<|a|<1$, we have that $C_{\omega, \psi}=C_{z^ne^{f(z)+g(\frac{1}{z})}, az}$ is isomorphically conjugate to $T_3^{-1}$ on $H(\mathbb{C}^*)$, where
\[
T_3=C_{cz^n, az}
\]
and $c=a^n e^{-g(0)}\ne 0$. Proposition~\ref{P:proposition17} and Remark~\ref{R:inverseFHC} give that both $T_3$ and $T_3^{-1}$ satisfy the Frequent Hypercyclicity Criterion, and hence so does $C_{\omega, \psi}$ by Remark~\ref{R:1.13}.

Similarly, if Condition (3)(ii) holds %, performing the isomorphic conjugacies as in the proof of the second case of $(1)\Rightarrow (3)$ 
we see that
 \[
C_{w,\psi}=C_{z^{-n} e^{f(z)+g(\frac{1}{z})}, \frac{1}{a} z}
\]
is isomorphically conjugate to $C_\frac{1}{z}C_{w,\psi}C_\frac{1}{z}=C_{z^ne^{f(\frac{1}{z})+g(z)},az}$, which by Lemma~\ref{L:corollary15} is isomorphically conjugate to
\[
T:= M_{z^{n} e^{f(\frac{1}{z})+g(0)}} C_{az}.
\]
Now, $T^{-1}$ is isomorphically conjugate to $C_\frac{1}{z}T^{-1}C_\frac{1}{z}=C_{a^n z^n e^{-f(az)-g(0)}, az}$, which by Lemma~\ref{L:corollary15} is isomorphically conjugate to
\[
S:= C_{\tilde{c} z^n, az},
\]
where $\tilde{c}=a^n e^{-g(0)}\ne 0$, and 
by Proposition~\ref{P:proposition17} and  Remark~\ref{R:inverseFHC}  both $S$ and $S^{-1}$ satisfy the Frequent Hypercyclicity Criterion. But  $C_{\omega, \psi}$ is isomorphically conjugate (to $T$, and hence) to $S^{-1}$, so it   must satisfy the Frequent Hypercyclicity Criterion as well. 
\end{proof}

\setcounter{section}{3} 

%
%
%
%.  SUBSECTION
%
%
%

\subsection{The punctured disc case. Proof of Theorem~\ref{T:mixing,dstar} }

\setcounter{section}{1} 

\begin{theorem} %\label{T:mixing,dstar} 
Let $\omega:\mathbb{D}^*\to\mathbb{C}$ and $\psi:\mathbb{D}^*\to\mathbb{D}^*$ be holomorphic, where $\mathbb{D}^*=\mathbb{D}\setminus\{ 0\}$. The following are equivalent:
\begin{enumerate}
\item[(1)]\ $C_{\omega, \psi}$ is hypercyclic on $H(\mathbb{D}^{*})$.

\item[(2)]\ $C_{\omega, \psi}$ satisfies the Frequent Hypercyclicity Criterion on $H(\mathbb{D}^{*})$, and in particular is frequently hypercyclic, mixing and chaotic.

\item[(3)] The multiplier $\omega$ is of the form $\omega (z)=z^n e^{W(z)}$ $(z\in \mathbb{D}^*)$, where $W\in H(\mathbb{D}^*)$ and $n\in\mathbb{N}$, and the symbol $\psi $ is injective, not surjective, and satisfies $\lim_{z\to 0}\psi(z)=0$.

\end{enumerate}
\end{theorem}

The equivalence $(1)\Leftrightarrow (3)$ is due to Goli\'nski and Przestacki ~\cite[Theorem~{5.11}]{golprz2021}.
We use the following result \cite[p.93]{ShapiroCompOp}.
\setcounter{section}{3} 
\setcounter{theorem}{0}
\begin{theorem} {\bf (K\"{o}nig)} \label{T:Konig}

Let $\psi:\mathbb{D}\to\mathbb{D}$ be a univalent and non-surjective holomorphic map with $\psi(0)=0$. Then there exists $a\in\mathbb{D}^*$ and $\sigma\in H(\mathbb{D})$ univalent satisfying $\sigma(0)=0$ so that the diagram
\begin{equation} \label{eq:Konig}
\begin{CD}
\mathbb{D}@> \psi >>\mathbb{D}  \\
@V \sigma VV @ VV \sigma V \\
\sigma(\mathbb{D}) @ > az >> \sigma(\mathbb{D}).
\end{CD}
%\begin{tikzcd}
%\mathbb{D} \arrow{r}{\psi} \arrow[swap]{d}{\sigma} & \mathbb{D} \arrow{d}{\sigma} \\%
%\sigma(\mathbb{D}) \arrow{r}{az}& \sigma(\mathbb{D})
%\end{tikzcd}
\end{equation}
commutes.
\end{theorem}

%\begin{proposition} \label{proposition21}{\bf \cite[Proposition~5.3]{golprz2021}}

% Let $\omega(z) = F(\frac{1}{z})$, where $F \in H(\mathbb{C})$. Then $C_{\omega, az}$ is not hypercyclic on $H(\mathbb{D}^{*})$.  
%\end{proposition}

\begin{proof}[Proof of Theorem \ref{T:mixing,dstar}] 
The implication $(2)\Rightarrow (1)$ is immediate, and  $(1)\Leftrightarrow (3)$ is \cite[Theorem~{5.11}]{golprz2021}, so we show
$(3)\Rightarrow (2)$. Suppose $(3)$ holds. Notice that the isolated singularity of $\psi$ at $0$ must be removable, and %since $\psi$ is non-constant by the Maximum Modulus Principle we may assume that $\psi\in H(\mathbb{D})$ and 
by $(3)$ we may assume $\psi$ is a non-surjective univalent self-map of $\mathbb{D}$ satisfying $\psi(0)=0$.
By K\"{o}nig's theorem there exist $\sigma:\mathbb{D}\to\mathbb{C}$ univalent with $\sigma(0)=0$ and $a\in\mathbb{D}^*$ satisfying $\eqref{eq:Konig}$ and thus the diagram 
\[
\begin{CD}
\mathbb{D}^* @>   \psi >> \mathbb{D}^*  \\
@V \sigma VV @ VV \sigma V \\
\Omega^* @ > az >> \Omega^*.
\end{CD}
\]
%\[
%\begin{tikzcd}
%\mathbb{D}^* \arrow{r}{\psi} \arrow[swap]{d}{\sigma} & \mathbb{D}^* \arrow{d}{\sigma} \\%
%\Omega^* \arrow{r}{az}& \Omega^*
%\end{tikzcd}
%\]
commutes, where $\Omega=\sigma(\mathbb{D})$, $\Omega^*=\Omega\setminus\{0\}$, and $a=\psi'(0)$. 
Hence the diagram
\[
\begin{CD}
H(\Omega^*)@> C_{\omega\circ\sigma^{-1}, az} >>H(\Omega^*) \\
@V C_\sigma VV @ VV C_\sigma V \\
H(\mathbb{D}^*) @ > C_{\omega, \psi} >> H(\mathbb{D}^*).
\end{CD}
\]
%
%\[
%\begin{tikzcd}
% H(\Omega^*) \arrow{r}{C_{\omega\circ\sigma^{-1} ,az}} \arrow[swap]{d}{C_\sigma} &  H(\Omega^*) \arrow{d}{C_\sigma} \\%
% H(\mathbb{D}^*) \arrow{r}{C_{\omega,\psi}} & H(\mathbb{D}^*) 
%\end{tikzcd}
%\]
commutes, with $C_\sigma$ invertible. So it suffices to show that $C_{\omega\circ \sigma^{-1}, az}$ satisfies the Frequent Hypercyclicity Criterion on $H(\Omega^*)$, by Remark~\ref{R:1.13}. Now, by $(3)$ there exist $n\in\mathbb{N}$ and $W\in H(\mathbb{D}^*)$ so that
\begin{equation}
\omega(\xi)=\xi^n  e^{W(\xi)} \ \ \ \xi \in \mathbb{D}^*.
\end{equation}
Also, since $\mathbb{D}$ is simply connected and $\sigma:\mathbb{D}\to \mathbb{C}$ univalent with $\sigma(0)=0$, there exists $g\in H(\mathbb{D})$ so that 
\[
\sigma(\xi)=\xi e^{g(\xi)} \ \ \ \ (\xi \in \mathbb{D}). 
\]
Now, since each $z\in \Omega^*$ is uniquely represented as $z=\sigma(\xi)$ with $\xi\in \mathbb{D}^*$ we may write
\[
\begin{aligned}
(\omega\circ \sigma^{-1})(z)&= (\omega\circ \sigma^{-1})(\sigma(\xi)) =\omega(\xi) \\
&= \xi^n e^{W(\xi)}\\
& =\left( \frac{\sigma(\xi)}{e^{g(\xi)}} \right)^n e^{W\circ \sigma^{-1}(z)}\\
&= z^n e^{h(z)}
\end{aligned}
\]
where $h=-n(g\circ \sigma^{-1})+W\circ\sigma^{-1}\in H(\Omega)$. But by Lemma~\ref{L:corollary15} the operator $C_{\omega\circ\sigma^{-1}, az}=C_{z^n e^{h(z)}, az}$ is isomorphically conjugate to \[ C_{z^ne^{h(0)}, az}:H(\Omega^*)\to H(\Omega^*)\] which satisfies the Frequent Hypercyclicity Criterion by Proposition~\ref{P:proposition17}. 
 \end{proof}

\begin{proposition}\label{P:noscDast}
Let $\psi:\mathbb{D}^*\to \mathbb{D}^*$ be holomorphic and surjective, and let $\omega\in H(\mathbb{D}^*)$. Then $C_{\omega, \psi}$ is not supercyclic on $H(\mathbb{D}^*)$.
\end{proposition}

\begin{proof}
Suppose $C_{\omega, \psi}$ is supercyclic on $H(\mathbb{D}^*)$. By Theorem~\ref{T:necessary}, $\omega$ is zero-free and $\psi:\mathbb{D}^*\to \mathbb{D}^*$ is univalent, strongly run-away and $\mathbb{D}^*$-convex. Arguing as in the proof of Theorem~\ref{T:mixing,dstar}, we may assume $\psi$ is a holomorphic self-map of $\mathbb{D}$. 

We claim that $\psi(0)=0$. Indeed, suppose that $b:=\psi(0)\ne 0$. Since $\psi$ is non-constant, there exists $0<\rho<1$ so that 
$
\psi(z)\ne b
$
for any $z$ with $0<|z|\le \rho$. By the Inverse Function Theorem \cite[p 234]{gamelin2001}, for 
\[
0<\delta < \mbox{min}\{ |\psi(z)-b|: \ z\in C(0, \rho)\}
\]
we know that $\psi^{-1}$ is univalent on $D(b,\delta)$ and satisfies $\psi^{-1} (D(b, \delta))\subset D(0, \rho)$ and $\psi^{-1}(b)=0$.
So 
\[
\gamma:=\psi^{-1} (C(b, \frac{\delta}{2}))
\]
is a compact jordan arc contained in $\mathbb{D}^*$ that is $\mathbb{D}^*$-convex. Indeed, since univalent maps preserve orientation and $\psi^{-1}(D(b, \frac{\delta}{2}))$ is the component of $\psi^{-1} (D(b,\delta)\setminus C(b, \frac{\delta}{2}))$ that contains $0=\psi^{-1}(b)$, the latter belongs to the hole of $\gamma$.  That is, the image of the $\mathbb{D}^*$-convex compact subset $\gamma$ of $\mathbb{D}^*$ is the compact set $\psi(\gamma)=C(b, \frac{\delta}{2})$ which is contained in $D(b, \delta)$ and thus  is not $\mathbb{D}^*$-convex. This contradicts that the map $\psi$ is $\mathbb{D}^*$-convex, so the claim that $\psi(0)=0$ holds. 
Hence we now have  by our surjectivity assumption that $\psi\in \mbox{Aut}(\mathbb{D})$ with $\psi(0)=0$, which forces $\psi$ to be a rotation around $0$ or the identity,
contradicting in either case that $\psi$ is run-away on $\mathbb{D}^*$.
\end{proof}

Theorem~\ref{T:supercyclic,cstar}, Proposition~\ref{P:noscDast}, Remark~\ref{R:1.13} and the Riemann mapping theorem now give the following.

\begin{corollary} \label{C:sc=FHCpunctured}
Let $\Omega\subseteq \mathbb{C}$ be a simply connected domain and $z_0\in \Omega$, and let $\Omega^*:=\Omega\setminus\{ z_0\}$. Then for any $\psi\in \mbox{Aut}(\Omega^*)$ and $\omega\in H(\Omega^*) $, the composition operator $C_{\omega, \psi}:H(\Omega^*)\to H(\Omega^*)$  is supercyclic if and only if it satisfies the Frequent Hypercyclicity Criterion.
\end{corollary}

%\section{Finitely connected domains conformally equivalent to an annulus or having at least two holes}
\section{The case when $\Omega$ is an annulus or has two or more holes} \label{S:2ormore}
\subsection{Proof of Theorem~\ref{T:NoSChere}}

According to Theorem~\ref{T:NoSChere}, when $\Omega$ is conformally equivalent to a non-degenerate annulus or is $n$-connected for some $n\ge 3$, the space $H(\Omega)$ does not support supercyclic weighted composition operators. 
This will follow once we establish Proposition~\ref{P:annulus,supercyclic} and Proposition~\ref{P:3ormore, supercylic}  below. 

\begin{proposition} \label{P:annulus,supercyclic} 
Let $A(r) = \{z \in \mathbb{C}: 1 < \lvert z \rvert < r\}$, where $r>1$, and let $\omega: A(r) \to \mathbb{C}$ and $\psi: A(r) \to A(r)$ be holomorphic. Then $C_{\omega, \psi}$ is not supercyclic on $H(A(r))$.
\end{proposition}

The proof of Proposition~\ref{P:annulus,supercyclic} relies on the following observation by Goli\'nsky and Przestacki. 

\begin{lemma} \label{L:annulus,automorphism} {\rm (\cite[Section 6]{golprz2021})} If $\psi: A(r) \to A(r)$ is univalent and the hole of A(r) is contained in the hole of $\psi(A(r))$, then $\psi$ is an $A(r)$-automorphism.
\end{lemma}

\begin{proof}[Proof of Proposition~\ref{P:annulus,supercyclic}] Suppose $C_{\omega, \psi}$ is supercyclic on $H(A(r))$. Then $\psi$ is an $A(r)$-automorphism, by
 Lemma \ref{L:supercyclic,necessary,complex}, Lemma \ref{L:supercyclic,convex}, and Lemma~\ref{L:annulus,automorphism}. So $\psi$ is the composition of a rotation and (optionally) the inversion map $z \to \frac{r}{z}$, forcing that for any $z \in A(r)$ the set $\overline{\{\psi_n (z): n \in \mathbb{N}\}}$ is compact in $A(r)$, contradicting Lemma~\ref{L:supercyclic,necessary,complex}.
\end{proof}

We next consider the case when $\Omega$ is $n$-connected for some $n\ge 3$.

\begin{proposition} \label{P:3ormore, supercylic} 
Let $\Omega \subset \mathbb{C}$ be a finitely connected domain with at least two holes. For all holomorphic functions $\omega: \Omega \to \mathbb{C}$ and $\psi: \Omega \to \Omega$, the composition operator $C_{\omega, \psi}$ is not supercyclic on $H(\Omega)$. 
\end{proposition}

The proof of Proposition~\ref{P:3ormore, supercylic} relies on the following.% from \cite{GolPrzCharacterization2021}.

\begin{lemma} \label{L:3ormore,automorphism} {\bf \cite[Section 6]{golprz2021}}
Let $\Omega \subset \mathbb{C}$ be a finitely connected domain with at least two holes. Then
\begin{enumerate}
\item[(i)]\ the group $\mbox{Aut}(\Omega)$ of automorphisms on $\Omega$ is finite, and
\item[(ii)]\ if $\psi: \Omega \to \Omega$ is univalent and each hole of $\Omega$ is contained in a hole of $\psi(\Omega)$, then $\psi\in \mbox{Aut}(\Omega)$. 
\end{enumerate}
\end{lemma}

\begin{proof}[Proof of Proposition~\ref{P:3ormore, supercylic}] Assume that $C_{\omega, \psi}$ is supercyclic on $H(\Omega)$. Then $\psi\in \mbox{Aut}(\Omega)$ and $\mbox{Aut}(\Omega)$ is finite, by Lemma \ref{L:supercyclic,necessary,complex}, Lemma \ref{L:supercyclic,convex}, and Lemma~\ref{L:3ormore,automorphism}, so for any fixed $z \in \Omega$, the set $\{\psi_n (z): n \in \mathbb{N}\}$ is finite, contradicting Lemma \ref{L:supercyclic,necessary,complex}.
\end{proof}

\subsection{Infinitely connected domains}

We show in this subsection Theorem~\ref{T:supercyclic=mixing,infinite}, that when $\Omega$ is infinitely connected the operator  $C_{\omega, \psi}$ is supercyclic on $H(\Omega)$ if and only if it is mixing. 

\setcounter{section}{1}
\setcounter{theorem}{5}
\begin{theorem} %\label{T:supercyclic=mixing,infinite}
 Let $\omega: \Omega \to \mathbb{C}$ and $\psi: \Omega \to \Omega$ be holomorphic, where $\Omega \subset \mathbb{C}$ is an infinitely connected domain. The following are equivalent: 
\begin{enumerate}
\item[(1)]\ $C_{\omega, \psi}$ is supercyclic on $H(\Omega)$. 

\item[(2)]\ $C_{\omega, \psi}$ is mixing on $H(\Omega)$. 

\item[(3)]\ The multiplier $\omega$ is zero-free and the symbol $\psi$ is univalent, $\Omega$-convex and strongly run-away.
\end{enumerate}
\end{theorem}
\setcounter{section}{4}

The proof of Theorem~\ref{T:supercyclic=mixing,infinite} relies on a lemma similar to \cite[Lemma~{3.13}]{GroMor2009}. 

\begin{lemma} \label{L:convexity,infinite} Let $\Omega \subset \mathbb{C}$ be a domain, and let $\psi: \Omega \to \Omega$ be univalent, $\Omega$-convex and strongly run-away.  Then for every connected $\Omega$-convex compact subset $K \subset \Omega$ with at least 2 holes, there exists $N_1 \in \mathbb{N}$ such that for each $n \geq N_1$ the set $\psi_n (K) \cup K$ is $\Omega$-convex. 
\end{lemma}

\begin{proof} Suppose otherwise. Then there exists an $\Omega$-convex connected compact subset $K$ of $\Omega$ and a  subsequence $(m_j)$ of $(n)$  such that for every $j \in \mathbb{N}$, $\psi_{m_j} (K) \cup K$ is not $\Omega$-convex. Fix an exhaustion $(K_i)$ of $\Omega$ of connected, $\Omega$-convex compact sets, all containing K. Note that since $\psi$ is $\Omega$-convex and univalent, each $\psi_{m_j} (K)$ is $\Omega$-convex and has the same number of holes as K has, which by assumption is at least 2. 

We have three possible cases:

\noindent Case 1: If for some $j \in \mathbb{N}$, $\psi_{m_j} (K)$ lies in the unbounded component of $\mathbb{C} \backslash K$ and K lies in the unbounded component of $\mathbb{C} \backslash \psi_{m_j} (K)$, then it follows that $\psi_{m_j} (K) \cup K$ is $\Omega$-convex, a contradiction. So this can't happen.
               %%%%%%%%%%%%%%%%%%%%%%%%%%%%%%%%%%%%%%%%%%%%AQUI!!!!!!!!!!!
               
\noindent Case 2: Infinitely many of the $\psi_{m_j} (K)$ lie in holes of K. By Lemma~\ref{L:cons}(i) the set K has only finitely many holes, so infinitely many of the $\psi_{m_j} (K)$ lie in some fixed hole $\mathcal{O}$ of K. By passing to a subsequence if necessary, we may assume that all of them do. Let $i \in \mathbb{N}$ such that $\psi_{m_1} (K) \subset K_i$. Since $\psi$ is strongly run-away we may choose $j \in \mathbb{N}$ large enough so that $\psi_{m_j} (K_i) \cap K_i = \emptyset$. Notice that $\psi_{m_1} (K)$ and $\psi_{m_j} (K)$ are disjoint subsets of $\mathcal{O}$. We have three possibilities. First, if both $\psi_{m_1} (K)$ and $\psi_{m_j} (K)$ lie in the unbounded component of the complement in $\mathbb{C}$ of the other, then both $\psi_{m_j} (K) \cup K$ and $\psi_{m_1} (K) \cup K$ are $\Omega$-convex, a contradiction. Second, if $\psi_{m_1} (K)$ lies in a hole of $\psi_{m_j} (K)$, then $\psi_{m_1} (K) \cup K$ is $\Omega$-convex because $\psi_{m_j} (K)$ has at least two holes, a contradiction. Third, if $\psi_{m_j} (K)$ lies in a hole of $\psi_{m_1} (K)$, then $\psi_{m_j} (K) \cup K$ is $\Omega$-convex because $\psi_{m_1} (K)$ has at least two holes, a contradiction. 

\noindent Case 3: For infinitely many $j \in \mathbb{N}$, K lies in holes of $\psi_{m_j} (K)$. Again we may assume this is true for all j. Let $i \in \mathbb{N}$ such that $\psi_{m_1} (K) \subset K_i$. Choose $j \in \mathbb{N}$ large enough so that $\psi_{m_j} (K_i) \cap K_i = \emptyset$. Notice that $\psi_{m_1} (K)$ and $\psi_{m_j} (K)$ are disjoint sets. Since both of these sets contain K in one of their holes, either $\psi_{m_1} (K)$ lies in a hole of $\psi_{m_j} (K)$ or $\psi_{m_j} (K)$ lies in a hole of $\psi_{m_1} (K)$. We then argue as in Case 2, obtaining a contradiction with both possibilities. 
\end{proof}

\begin{proof}[Proof of Theorem~\ref{T:supercyclic=mixing,infinite}]  The implication $(2)\Rightarrow (1)$ is immediate, and $(1)\Rightarrow (3)$ follows by Theorem~\ref{T:necessary}. To see  $(3)\Rightarrow (2)$, let U and V be non-empty open subsets of $H(\Omega)$. From the topology on $H(\Omega)$, there exists $\epsilon > 0$, a compact subset K of $\Omega$, and functions $f, g \in H(\Omega)$ such that
\[
\{h \in H(\Omega): \lvert \lvert h - f \rvert \rvert_K < \epsilon \} \subset U 
\]
\noindent and 
\[
\{h \in H(\Omega): \lvert \lvert h - g \rvert \rvert_K < \epsilon \} \subset V.
\]
Enlarging K if necessary, we may assume K is a connected $\Omega$-convex compact subset of $\Omega$ with at least two holes. So by Lemma \ref{L:convexity,infinite}, there exists $N \in \mathbb{N}$ such that for each $n \geq N$, 
the set $\psi_n (K) \cup K$ is $\Omega$-convex and $\psi_n (K) \cap K = \emptyset$. By Runge's theorem there exists $h_n \in H(\Omega)$ such that 
\[
\lvert \lvert h_n (z) - f(z) \rvert \rvert_K < \epsilon 
\]
\noindent and 
\[
\lvert \lvert h_n (z) - \frac{g (\psi_{-n} (z))}{\prod_{i=1}^{n} \omega (\psi_{-i} (z))} \rvert \rvert_{\psi_n (K)} < \frac{\epsilon}{\lvert \lvert  \prod_{i=0}^{n-1} \omega (\psi_i (z)) \rvert \rvert_K}. 
\]

\noindent So $h_n \in U$,  and from the second inequality, we get
\begin{flalign*}
\left \| \frac{(\prod_{i=1}^{n} \omega (\psi_{-i} (z))) h_n (z) - g (\psi_{-n} (z))}{\prod_{i=1}^{n} \omega (\psi_{-i} (z))} \right \|_{\psi_n (K)} 
& = \left\|  \frac{(\prod_{i=0}^{n-1} \omega (\psi_i (z))) h_n (\psi_n (z)) - g (z)}{\prod_{i=0}^{n-1} \omega (\psi_i (z))} \right\|_K&\\ 
& = \left \| \frac{(C_{\omega, \psi}^{n} (h_n)) (z) - g (z)}{\prod_{i=0}^{n-1} \omega (\psi_i (z))} \right \|_K &\\ 
& < \frac{\epsilon}{\lvert \lvert \prod_{i=0}^{n-1} \omega (\psi_i (z)) \rvert \rvert_K}.
\end{flalign*}
\noindent So $C_{\omega, \psi}^{n} (h_n) \in V$, and  $C_{\omega, \psi}$ is mixing.
\end{proof}

\begin{proof}[Proof of Theorem~\ref{T:overall}]
The conclusion is known for the case when $\Omega$ is simply connected \cite[Theorem 2.1]{BesCAOT}. The punctured simply connected case follows from Theorem~\ref{T:supercyclic,cstar}, Theorem~\ref{T:mixing,dstar}, and Corollary~\ref{C:sc=FHCpunctured}. The remaining cases follow from Theorem~\ref{T:NoSChere} and Theorem~\ref{T:supercyclic=mixing,infinite}.
\end{proof}

\begin{remark}
{\rm We were not able to show the equivalence of supercyclicity and mixing for 
$C_{\omega, az}: H(\mathbb{D}^{*}) \to H(\mathbb{D}^{*})$. This equivalence between supercyclicity and mixing would hold if and only if the following conjecture is true.}
\begin{quote} {\bf Conjecture.} \ Let $\Omega$ be a simply connected domain containing zero and $a\in\mathbb{D}^*$ so that $a\Omega \subseteq \Omega$.
For each $F \in H(\mathbb{C})$, the operator $C_{F(\frac{1}{z}), az}$ is not supercyclic  on $H(\Omega^*)$.
\end{quote}
\end{remark}
Indeed, it suffices to show the conjecture for the case when $F$ is not a polynomial. A positive answer to the conjecture ensures, together with the results from this paper, that for any planar domain $\Omega$,  a weighted composition operator on $H(\Omega)$ is supercyclic if and only if it is mixing.

%\section{Some open problems}

%\newpage

\end{document}